\documentclass{article}

\usepackage{arxiv}

\usepackage[utf8]{inputenc} 
\usepackage[T1]{fontenc}    
\usepackage{hyperref}       
\usepackage{url}            
\usepackage{booktabs}       
\usepackage{amsfonts}       
\usepackage{nicefrac}       
\usepackage{microtype}      
\usepackage{lipsum}
\usepackage{amsmath}

 \newtheorem{thm}{Theorem}[section]

 \newtheorem{defn}[thm]{Definition}
 \newtheorem{rem}[thm]{Remark}

 \numberwithin{equation}{section}
 
\title{On The Horadam Hybrid Quaternions}

\author{
  Ali Dağdeviren\\
  Weight and Balance Department\\
  Turkish Aviation Academy\\
  Istanbul, Bakırköy, 34149 \\
  \texttt{adagdeviren@thy.com} \\
   \And
 Ferhat Kürüz \\
  Department of Computer Engineering\\
  Faculty of Engineering and Architecture\\
  Istanbul Gelisim University \\
  Istanbul, Bakırköy, 34149\\
  \texttt{fkuruz@gelisim.edu.tr} \\
}

\begin{document}
\maketitle

\begin{abstract}
In this paper, we define the Horadam hybrid quaternions and give some of their properties. Moreover, we investigate the relations
between the Fibonacci hybrid quaternions and the Lucas hybrid  quaternions which connected
the Fibonacci quternions and Lucas quaternions. Furthermore, we give the Binet formulas and Cassini identities for these quaternions.
\end{abstract}


\section{Introduction}
Quaternions were introduced by Sir William Rowan Hamilton in 1866 as an extension of the complex numbers \cite{quat}. Quaternions are an important number system used in different areas such as computer science, quantum physics, and analysis \cite{quatcomp2, quatcomp4,quatcomp1}. This type of quaternions also called real quaternions. A real quaternion is defined as
\begin{equation*}
    Q=z_{0}+z_{1}i+z_{2}j+z_{3}k
\end{equation*}
where $z_0, z_1, z_2, z_3$ are real numbers. Also, $i$, $j$, and $k$ are the units of real quaternions which satisfy the following equalities
\begin{equation}\label{1}
    i^2=j^2=k^2=ijk=-1.
\end{equation}
Note that the set of all quaternions form an associative but non-commutative algebra. The conjugate of the quaternion $Q$ is denoted by  $\overline{Q}$ and defined by $\overline{Q}=z_{0}-z_{1}i-z_{2}j-z_{3}k$. Moreover, the norm of any quaternion $Q$ is denoted by $\left\| Q \right\|$ and defined by $\left\| Q \right\|=\sqrt{Q\overline{Q}}=\sqrt{z_0^2+z_1^2+z_2^2+z_3^2}\in \mathbb{R}$. For more information about quaternions see \cite{quat2}.

Hybrid numbers have been defined by {\"O}zdemir \cite{hybnum} as a new generalization of complex, dual and hyperbolic numbers. In this generalization, the author gave a number system of such numbers that consist of all three number systems together. The set of hybrid numbers is denoted by $\mathbb{K}$ and defined as follows:
 \begin{equation}\label{4}
    \mathbb{K}=\bigg\{ z=a+b\textbf{i}+c\boldsymbol{\varepsilon}+d\textbf{h}: ~a,b,c,d \in \mathbb{R},
    \begin{array}{l@{}c@{}}
     &\textbf{i}^2=-1, \boldsymbol{\varepsilon}^2=0, \textbf{h}^2=1 \\
     & \textbf{i}\textbf{h}=-\textbf{h}\textbf{i}=\boldsymbol{\varepsilon}+\textbf{i}
     \end{array}
    \bigg\}.
\end{equation}
Multiplication rules of  $ \textbf{i}, \boldsymbol{\varepsilon}$,  and $ \textbf{h}$ can be given as following table:
\smallskip

\begin{center}
\begin{tabular}{|c||c c c|}
\hline

\hline
$ . $ & \textbf{i} & $\boldsymbol{\boldsymbol{\varepsilon}}$ & $\textbf{h}$ \\

\hline \hline
 $\textbf{i}$ & $-1$  & $1-\textbf{h}$ & $\boldsymbol{\boldsymbol{\varepsilon}}$ + $\textbf{i}$ \\
 $\boldsymbol{\boldsymbol{\varepsilon}}$ & $1+\textbf{h}$& $0$ & $-\boldsymbol{\boldsymbol{\varepsilon}}$ \\
 $\textbf{h}$ & $-\boldsymbol{\boldsymbol{\varepsilon}} -\textbf{i}$  & $\boldsymbol{\boldsymbol{\varepsilon}}$  & $1$ \\

\hline

\hline
\end{tabular}
\smallskip\\
\small Table 1. Multiplication table of hybrid units
\end{center}
The conjugate of the hybrid number $z$ is denoted by ${z}^c$ and defined by ${z}^c=a - b \textbf{i} - c \boldsymbol{\varepsilon} -d \textbf{h}$. The norm of any hybrid number $z$ is
\begin{equation*}
    \left\|z\right\|=\sqrt{z{z}^c}=\sqrt{a^2+(b-c)^2-c^2-d^2}.
\end{equation*}
For further information about hybrid number system, see \cite{hybnum}.

The hybrid quaternions have recently been defined as a new quaternion system in \cite{Dagdeviren}. This system has a strong algebraic structure and it is a generalization of complex, dual, and hyperbolic quaternions. That is, complex quaternions, dual quaternions and hyperbolic quaternions can be obtained from hybrid quaternions in special cases. Moreover, hybrid quaternions are also generalised the features of the other three quaternion systems such as inner product, vector product, and norm. The set of hybrid quaternions denoted by $\mathbf{H}_\mathbb{K}$ and defined as
\begin{equation}
    \mathbf{H}_\mathbb{K}=\{Q=z_0+z_1 i+z_2 j+z_3 k ~|~ z_0, z_1, z_2, z_3 \in \mathbb{K}\}
\end{equation}
where quaternionic units $i,j,k$ satisfies the equation \eqref{1}. The quaternionic units $i,j,k$ commutes with the hybrid units $\textbf{i},\boldsymbol{\varepsilon}, \textbf{h}$. Thus, any hybrid quaternion can be written as
\begin{equation*}
    Q= q_0+ q_1 \textbf{i} + q_2 \boldsymbol{\varepsilon} + q_3 \textbf{h}
\end{equation*}
where $q_0, q_1, q_2, q_3$ are quaternions and $\textbf{i}, \boldsymbol{\varepsilon}, \textbf{h}$ are hybrid units obeying the multiplication rules in the Table 1, \cite{Dagdeviren}.

Two hybrid quaternions are equal if all their components are equal, one by one. Let $Q=q_0+q_1 \textbf{i}+q_2 \boldsymbol{\varepsilon} +q_3\textbf{h}$ and $P=p_0+ p_1 \textbf{i}+ p_2 \boldsymbol{\varepsilon} +p_3 \textbf{h}$ and  be any two hybrid quaternions. Addition and subtraction of these two hybrid quaternions are defined as
\begin{equation*}
    Q \mp P = (q_0 \mp p_0) +  (q_1 \mp p_1)\textbf{i} + (q_2 \mp p_2) \boldsymbol{\varepsilon} + (q_3 \mp p_3) \textbf{h}.
\end{equation*}
Multiplication of the hybrid quaternions is defined as
\begin{align*}
    QP &= (q_0 p_0-q_1p_1+q_3p_3+q_1p_2+q_2p_1)\\
      & \quad +(q_0p_1+q_1p_0+q_1p_3-q_3p_1) ~ \textbf{i} \\ 
      & \quad +(q_0p_2+q_2p_0+q_1p_3-q_3p_1+q_3p_2-q_2p_3)   ~\boldsymbol{\varepsilon}\\
      & \quad +(q_0p_3+q_3p_0+q_2p_1-q_1p_2)~ \textbf{h}~.
\end{align*}
\indent The product of any two hybrid numbers with the help of hybrid units is expressed as above. Furthermore, the hybrid quaternions $Q$ and $P$ above can be written as $Q=z_0+z_1i+z_2j+z_3k$ and $P=t_0+t_1i+t_2j+t_3k$ in terms of quaternion units. Multiplication of these two hybrid quaternions can be given as  
\begin{align*}
    QP &=       (z_0t_0-z_1t_1-z_2t_2-z_3t_3)\\
   & \quad +   (z_1t_0+z_0t_1-z_3t_2+z_2t_3)~ i\\ 
   & \quad +   (z_2t_0+z_3t_1+z_0t_2-z_1t_3)~ j\\
   & \quad +   (z_3t_0-z_2t_1+z_1t_2+z_0t_3)~ k~ .
\end{align*}
Further information about hybrid quaternions can be found in \cite{Dagdeviren}.

Generalized Fibonacci numbers were defined by Horadam \cite{horfib,hor}. Moreover Horadam gave the formula for negative index terms of Horadam numbers \cite{hor}. Horadam numbers are defined as follows:
\begin{equation}\label{horadam1}
    w_n=pw_{n-1}-qw_{n-2} ~ ,~ n\geq 2
\end{equation}
where $p$, $q$, $n$ are integers and $w_0$, $w_1$ are initial conditions. The Binet's formula of the Horadam numbers is
\begin{equation}
    w_n=A\alpha^n+B\beta^n
\end{equation}
where  $\alpha$ and $\beta$ are the roots of the equation $x^2-px+q=0$. Also $A$ and $B$ are 
\begin{equation}\label{horadam2}
    A=\frac{w_1-w_0\beta}{\alpha-\beta}~ ~ \text{,} ~~ B=\frac{w_0\alpha-w_1}{\alpha-\beta}~~~~.
\end{equation}
In addition, Horadam numbers can be represent as $w_n(w_0,w_1;p,q)$. For special $w_0,w_1,p$ and $q$ the equation \eqref{horadam1} defines the well known numbers named as the Fibonacci type numbers. These types of numbers can be listed as follows:
\begin{itemize}
    \item [i)] $w_n(0,1;p ,q) ~=~U_n-$ Generalized Fibonacci numbers,
    \item [ii)] $w_n(2,1;p ,q) ~=~V_n-$ Generalized Lucas numbers,
    \item [iii)] $w_n(0,1;1,-1)=F_n-$ Fibonacci numbers,
    \item [iv)] $w_n(2,1;1,-1)=L_n-$  Lucas numbers,
    \item [v)] $w_n(0,1;2,-1)=P_n-$  Pell numbers,
    \item [vi)] $w_n(2,2;2,-1)=PL_n-$  Pell-Lucas numbers,
    \item [vii)] $w_n(0,1;1,-2)=J_n~-$ Jacobsthal numbers,
    \item [viii)] $w_n(2,1;1,-2)=JL_n-$  Jacobsthal-Lucas numbers,
    \item [ix)] $w_n(0,1;3,~2) ~~=~M_n-$  Mersenne Numbers,
    \item [x)] $w_n(1,3;3,-2)=FE_n-$  Fermat numbers.
\end{itemize}

These numbers have been studied from different perspectives \cite{pelllucquat, horhyb, gendualpellquat, compfibquat, fibquat, hor,  pellpadquat, genkhor}. A special kind of hybrid numbers, namely Horadam hybrid numbers, were introduced and studied by Szynal-Liana \cite{horhyb}. In this study,  the author gave the Binet formulas, generating functions and characters for Horadam hybrid numbers. $n$th Horadam hybrid number is defined as
\begin{equation}\label{Horadam5}
    \breve{H}_n=w_n+ \textbf{i}~ w_{n+1} +\boldsymbol{\varepsilon} ~w_{n+2}+\textbf{h} ~w_{n+3}
\end{equation}
where $w_n$ is $n$th Horadam number. The Binet formula for the Horadam hybrid numbers is
\begin{equation}\label{22}
    w_n=A\alpha^*\alpha^n+B\beta^*\beta^n
\end{equation}
where $A, B$ are defined by \eqref{horadam2} and $\alpha^*=1+\mathbf{i}\alpha+\boldsymbol{\varepsilon}\alpha^2+\textbf{h}\alpha^3$, $\beta^*=1+\mathbf{i}\beta+\boldsymbol{\varepsilon}\beta^2+\textbf{h}\beta^3$. Moreover, in the same article, the author defined Fibonacci and Lucas hybrid numbers as follows:
\begin{align*}
    \breve{F}_n&={F}_n+\textbf{i}{F}_{n+1}+\boldsymbol{\varepsilon}{F}_{n+2}+\textbf{h}{F}_{n+3} ,\\
    \breve{L}_n&={L}_n+\textbf{i}{L}_{n+1}+\boldsymbol{\varepsilon}{L}_{n+2}+\textbf{h}{L}_{n+3} .
\end{align*}
The other studies about Horadam hybrid numbers can be found in \cite{comphor,genhybfib,pellhybrinom, jacluchyb,  fibhyb}.

In the literature, there are also studies on fibonacci quaternions. Horadam defined $n$th Fibonacci and Lucas quaternions as follows:
\begin{align*}
    \tilde{F}_n&=F_n+iF_{n+1}+jF_{n+2}+kF_{n+3} ,\\
    \tilde{L}_n&=L_n+iL_{n+1}+jL_{n+2}+kL_{n+3}
\end{align*}
where $i$, $j$ and $k$ are quaternion units which satisfy equations \eqref{1}. Moreover, $F_n$ and $L_n$ are the $n$th Fibonacci and Lucas numbers, respectively. In \cite{fibquat}, Halıcı have studied on Fibonacci quaternions and gave the generating functions and Binet formulas for Fibonacci and Lucas quaternions. For the other studies about Fibonacci and Lucas quaternions see \cite{splitfib, fibquat}.

In this study, we firstly define the Horadam hybrid quaternions. We then will give Binet formulas for these numbers. Furthermore we will examine the Fibonacci and Lucas hybrid quaternions in detail and consequently we will give some properties and identities of these numbers.

In what follows, to avoid confusion we use notation as properly as we can. That's why the hat sign is going to be used for the all type of hybrid quaternions such as Fibonacci hybrid quaternions$( \widehat{F}_n )$ or Lucas hybrid quaternions $( \widehat{L}_n )$ and the tilde sign is going to be used for the other type quaternion numbers such as Fibonacci quaternions $( \tilde{F}_n )$ or Lucas quaternions $( \tilde{L}_n )$. Moreover, Fibonacci and Lucas hybrid numbers will be shown as $( \breve{F}_n )$,  $( \breve{L}_n )$ respectively. In addition to the explanation above, we will give the following table to make the symbols used in this study easier to understand.

\begin{center}
\begin{tabular}{c|c}
\hline
\textbf{Notations}  & \textbf{Numbers} \\

\hline 
 \\[-8pt] $\widehat{F}_n$ & Fibonacci hybrid quaternions \\
 $\widehat{L}_n$ & Lucas hybrid quaternions  \\
\hline
\\[-8pt]
 $\tilde{F}_n$ & Fibonacci quaternions    \\
 $\tilde{L}_n$ & Lucas quaternions     \\
\hline
\\[-8pt]
 $\breve{F}_n$ & Fibonacci hybrid numbers \\
 $\breve{L}_n$ & Lucas hybrid numbers   \\
\hline
\\[-8pt]
 $F_n$ & Fibonacci numbers \\
 $L_n$ & Lucas numbers   \\
\hline
\end{tabular}
\smallskip\\
\small \textbf{Table 2.} Notation table of Numbers
\end{center}

\section{Horadam Hybrid Quaternions}

In this section, we define Horadam hybrid quaternions by using Horadam numbers. Therefore, we define the Horadam hybrid quaternion $\widehat{H}_n$ as
\begin{equation}\label{21}
    \widehat{H}_n=\breve{H}_n+i\breve{H}_{n+1}+j\breve{H}_{n+2}+k\breve{H}_{n+3}
\end{equation}
where $i$, $j$, $k$ are quaternion units which satisfy equations \eqref{1} and $H_n$ is the $nth$ Horadam hybrid number as in \eqref{Horadam5}. Moreover, every Horadam hybrid quaternion $\widehat{H}_n$ can be written as 
\begin{align*}
    \widehat{H}_n &=  (w_n+ \textbf{i} w_{n+1} +\boldsymbol{\varepsilon} w_{n+2}+\textbf{h} w_{n+3})\\
         & \quad +  (w_{n+1}+ \textbf{i} w_{n+2} +\boldsymbol{\varepsilon} w_{n+3}+\textbf{h} w_{n+4})~ i\\ 
         & \quad +  (w_{n+2}+ \textbf{i} w_{n+3} +\boldsymbol{\varepsilon} w_{n+4}+\textbf{h} w_{n+5})~ j\\
         & \quad +  (w_{n+3}+ \textbf{i} w_{n+4} +\boldsymbol{\varepsilon} w_{n+5}+\textbf{h} w_{n+6})~ k~ \\
         & =\tilde{H}_n+\textbf{i}\tilde{H}_{n+1}+\boldsymbol{\varepsilon} \tilde{H}_{n+2}+\textbf{h}\tilde{H}_{n+3}
\end{align*}
where $\textbf{i}, \boldsymbol{\varepsilon}, \textbf{h}$ are hybrid units and $\tilde{H}_n=w_n+ i w_{n+1} +jw_{n+2}+kw_{n+3}$ is the $n$th Horadam quaternion. The Fibonacci, Lucas, Pell, Jacobsthal, Pell-Lucas, and Jacobsthal-Lucas hybrid quaternions can be defined both by using Fibonacci$(\breve{F}_n)$, Lucas$(\breve{L}_n)$, Pell$(\breve{P}_n)$, Jacobsthal$(\breve{J}_n)$, Pell-Lucas$(\breve{PL}_n)$, and Jacobsthal-Lucas$(\breve{JL}_n)$ hybrid number coefficients  and by using \linebreak Fibonacci$(\tilde{F}_n)$, Lucas$(\tilde{L}_n)$, Pell$(\tilde{P}_n)$, Jacobsthal$(\tilde{J}_n)$, Pell-Lucas$(\tilde{PL}_n)$, and Jacobsthal-Lucas$(\tilde{JL}_n)$ quaternion coefficients respectively. We will define these hybrid quaternions as follows:

\begin{itemize}
   
    \item [\textbf{i)}] $nth$ Fibonacci hybrid quaternion $\widehat{F}_n$ is
    \begin{align*}
        \widehat{F}_n & =\breve{F}_n+i\breve{F}_{n+1}+j\breve{F}_{n+2}+k\breve{F}_{n+3} \\
        & =\tilde{F}_n+\textbf{i}\tilde{F}_{n+1}+\boldsymbol{\varepsilon} \tilde{F}_{n+2}+\textbf{h}\tilde{F}_{n+3},
    \end{align*}

    \item[\textbf{ii)}] $nth$ Lucas hybrid quaternion $\widehat{L}_n$ is
    \begin{align*}
        \widehat{L}_n & =\breve{L}_n+i\breve{L}_{n+1}+j\breve{L}_{n+2}+k\breve{L}_{n+3} \\
        & =\tilde{L}_n+\textbf{i}\tilde{L}_{n+1}+\boldsymbol{\varepsilon} \tilde{L}_{n+2}+\textbf{h}\tilde{L}_{n+3},
    \end{align*}
    
    \item[\textbf{iii)}] $nth$ Pell hybrid quaternion $\widehat{P}_n$ is
    \begin{align*}
        \widehat{P}_n & =\breve{P}_n+i\breve{P}_{n+1}+j\breve{P}_{n+2}+k\breve{P}_{n+3} \\
        & =\tilde{P}_n+\textbf{i}\tilde{P}_{n+1}+\boldsymbol{\varepsilon} \tilde{P}_{n+2}+\textbf{h}\tilde{P}_{n+3},
    \end{align*}
    
    \item[\textbf{iv)}] $nth$ Jacobsthal hybrid quaternion $\widehat{J}_n$ is
    \begin{align*}
        \widehat{J}_n & =\breve{J}_n+i\breve{J}_{n+1}+j\breve{J}_{n+2}+k\breve{J}_{n+3} \\
        & =\tilde{J}_n+\textbf{i}\tilde{J}_{n+1}+\boldsymbol{\varepsilon} \tilde{J}_{n+2}+\textbf{h}\tilde{J}_{n+3},
    \end{align*}
    
    \item[\textbf{v)}] $nth$ Pell-Lucas hybrid quaternion $\widehat{PL}_n$ is
    \begin{align*}
        \widehat{PL}_n & =\breve{PL}_n+i\breve{PL}_{n+1}+j\breve{PL}_{n+2}+k\breve{PL}_{n+3} \\
        & =\tilde{PL}_n+\textbf{i}\tilde{PL}_{n+1}+\boldsymbol{\varepsilon} \tilde{PL}_{n+2}+\textbf{h}\tilde{PL}_{n+3},
    \end{align*}
    
    \item[\textbf{vi)}] $nth$ Jacobsthal-Lucas hybrid quaternion $\widehat{JL}_n$ is
    \begin{align*}
        \widehat{JL}_n & =\breve{JL}_n+i\breve{JL}_{n+1}+j\breve{JL}_{n+2}+k\breve{JL}_{n+3} \\
        & =\tilde{JL}_n+\textbf{i}\tilde{JL}_{n+1}+\boldsymbol{\varepsilon} \tilde{JL}_{n+2}+\textbf{h}\tilde{JL}_{n+3}.
    \end{align*}
    
\end{itemize}

\begin{thm}
Let $n \in \mathbb{N}$, then the Binet formula for the Horadam hybrid quaternions is
\begin{equation*}
    \widehat{H}_n=A \alpha^* \underline{\alpha} \alpha^n + B \beta^* \underline{\beta} \beta^n 
\end{equation*}
where $A, B$ are defined by \eqref{horadam2} and 
\begin{align*}
    \alpha^*=1+\textbf{i}\alpha+\boldsymbol{\varepsilon}\alpha^2+h\alpha^3 \quad,&\quad \beta^*=1+\textbf{i}\beta+\boldsymbol{\varepsilon}\beta^2+h\beta^3,\\ \underline{\alpha}=1+i\alpha+j\alpha^2+k\alpha^3 \quad,&\quad \underline{\beta}=1+i\beta+j\beta^2+k\beta^3.
\end{align*}
\end{thm}

\textbf{\textit{Proof.}}
By using the definition of Horadam hybrid quaternions \eqref{21} and the Binet formula for the Horadam hybrid numbers \eqref{22}, we get
\begin{align*}
   \widehat{H}_n &= (A\alpha^n\alpha^*+B\beta^n\beta^*) + i (A\alpha^{n+1}\alpha^*+B\beta^{n+1}\beta^*) \\ 
 & \quad + j (A\alpha^{n+2}\alpha^*+B\beta^{n+2}\beta^*) + k (A\alpha^{n+3}\alpha^*+B\beta^{n+3}\beta^*) \\
 & = (A \alpha^n\alpha^*)(1+i\alpha+j\alpha^2+k\alpha^3) +(B \beta^n\beta^*)(1+i\beta+j\beta^2+k\beta^3) \\
 & = A \alpha^* \underline{\alpha} \alpha^n + B \beta^* \underline{\beta} \beta^n.
\end{align*}

\section{Fibonacci and Lucas Hybrid Quaternions}

In this section, we examine the Fibonacci and Lucas Hybrid quaternions in detail and give some properties. 
\begin{defn}
We denote the set of Fibonacci hybrid quaternios by ${FHQ}$ and define as follows:
 \begin{equation}
    {FHQ}=\{ \widehat{F}_n=\breve{F}_n+i\breve{F}_{n+1}+j\breve{F}_{n+2}+k\breve{F}_{n+3} ~|~ ~  \breve{F}_{n},  n \textrm{th} ~  \textrm { Fibonacci hybrid number} \}
\end{equation}
where $i,j,k$ are quaternionic units. Moreover, here  $n$th, $(n+1)th$, $(n+2)th$ and $(n+3)th$  Fibonacci hybrid numbers are
\begin{align}
    \breve{F}_n~&=F_n+\textbf{i}F_{n+1}+\boldsymbol{\varepsilon} F_{n+2}+\textbf{h}F_{n+3} , \label{12} \\
     \breve{F}_{n+1}&=F_{n+1}+\textbf{i}F_{n+2}+\boldsymbol{\varepsilon} F_{n+3}+\textbf{h}F_{n+4},\label{13}\\
     \breve{F}_{n+2}&=F_{n+2}+\textbf{i}F_{n+3}+\boldsymbol{\varepsilon} F_{n+4}+\textbf{h}F_{n+5},\label{14}\\
     \breve{F}_{n+3}&=F_{n+3}+\textbf{i}F_{n+4}+\boldsymbol{\varepsilon} F_{n+5}+\textbf{h}F_{n+6}\label{15}
\end{align}
where  $\textbf{i}, \boldsymbol{\varepsilon},$ and $\textbf{h}$ are the hybrid units. We will restate $\widehat{F}_n$ by using the equations \eqref{12}, \eqref{13}, \eqref{14} and \eqref{15} as below.
\begin{align}
      \widehat{F}_n&=(F_n+\textbf{i}F_{n+1}+ \boldsymbol{\varepsilon} F_{n+2}+\textbf{h}F_{n+3}) \nonumber \\
            &\quad +i(F_{n+1}+\textbf{i}F_{n+2}+ \boldsymbol{\varepsilon} F_{n+3}+\textbf{h}F_{n+4}) \nonumber\\
            &\quad +j(F_{n+2}+\textbf{i}F_{n+3}+ \boldsymbol{\varepsilon} F_{n+4}+\textbf{h}F_{n+5}) \nonumber \\
            &\quad +k(F_{n+3}+\textbf{i}F_{n+4}+ \boldsymbol{\varepsilon} F_{n+5}+\textbf{h}F_{n+6}) \nonumber \\ 
       \widehat{F}_n & =\tilde{F}_n+\textbf{i} \tilde{F}_{n+1}+ \boldsymbol{\varepsilon} \tilde{F}_{n+2}+\textbf{h} \tilde{F}_{n+3}
\end{align}
where $\tilde{F}_n={F}_n+i{F}_{n+1}+j{F}_{n+2}+k{F}_{n+3}$ is a  Fibonacci quaternion. Therefore, the set of Fibonacci hybrid quaternios ${FHQ}$  can be redefined  as
\begin{equation}
    {FHQ}=\bigg\{ \widehat{F}_n=\tilde{F}_n+\textbf{i} \tilde{F}_{n+1}+ \boldsymbol{\varepsilon} \tilde{F}_{n+2}+\textbf{h} \tilde{F}_{n+3}~ |
        \begin{array}{l@{}c@{}}
             &\tilde{F}_{n}, n\textrm{th} ~ \textrm { Fibonacci}\\
             &  \textrm {quaternion}
        \end{array}
    \bigg\}.
\end{equation}
\end{defn}
\begin{rem}
Every Fibonacci hybrid quaternion
\begin{equation*}
    \widehat{F}_n=\breve{F}_n+i\breve{F}_{n+1}+j\breve{F}_{n+2}+k\breve{F}_{n+3}
\end{equation*}
can be written as 
\begin{equation*}
    \widehat{F}_n =\tilde{F}_n+\textbf{i} \tilde{F}_{n+1}+ \boldsymbol{\varepsilon} \tilde{F}_{n+2}+\textbf{h} \tilde{F}_{n+3}.
\end{equation*}
\end{rem}

\begin{defn}
We denote the set of Lucas hybrid quaternios by ${LHQ}$ and define as
  \begin{equation}
    {LHQ}=\bigg\{ \widehat{L}_n=\breve{L}_n+i\breve{L}_{n+1}+j\breve{L}_{n+2}+k\breve{L}_{n+3} ~|
        \begin{array}{l@{}c@{}}
             &\breve{L}_{n}, n \textrm{th} ~ \textrm {Lucas}\\
             &  \textrm {hybrid number}
        \end{array}
    \bigg\}
\end{equation}
where $i,j,k$ are quaternionic units. Moreover, here $\breve{L}_n$ is $n$th Lucas hybrid number 
\begin{equation*}
    \breve{L}_n=L_n+\textbf{i}L_{n+1}+\boldsymbol{\varepsilon} L_{n+2}+\textbf{h}L_{n+3}
\end{equation*}

As with Fibonacci hybrid quaternions above, Lucas hybrid quaternions can be redefined in terms of Lucas quaternions by
\begin{equation}
    {LHQ}=\bigg\{ \widehat{L}_n=\tilde{L}_n+\textbf{i} \tilde{L}_{n+1}+ \boldsymbol{\varepsilon} \tilde{L}_{n+2}+\textbf{h} \tilde{L}_{n+3} |
        \begin{array}{l@{}c@{}}
             &\tilde{L}_{n}, n\textrm{th} ~ \textrm {Lucas}\\
             &  \textrm {quaternion}
        \end{array}
    \bigg\}
\end{equation}
where $\textbf{i},\boldsymbol{\varepsilon},\textbf{h}$ are hybrid units. Moreover, here $\tilde{L}_n$ is $n$th Lucas quaternion 
\begin{equation*}
    \tilde{L}_n=L_n+iL_{n+1}+ jL_{n+2}+kL_{n+3}.
\end{equation*}
\end{defn}
\smallskip
\begin{rem}
Every Lucas hybrid quaternion
\begin{equation*}
    \widehat{L}_n=\breve{L}_n+i\breve{L}_{n+1}+j\breve{L}_{n+2}+k\breve{L}_{n+3}
\end{equation*}
can be written as 
\begin{equation*}
    \widehat{L}_n=\tilde{L}_n+\textbf{i} \tilde{L}_{n+1}+ \boldsymbol{\varepsilon} \tilde{L}_{n+2}+\textbf{h} \tilde{L}_{n+3}.
\end{equation*}
\end{rem}
\smallskip
\begin{defn}
Let $\widehat{Q}_n$ and $\widehat{P}_n$ be $nth$ terms of the Fibonacci hybrid quaternion sequences such that 
\begin{align}
     \widehat{Q}_n&=\breve{Q}_n+i\breve{Q}_{n+1}+j\breve{Q}_{n+2}+k\breve{Q}_{n+3}\\
   & = \tilde{Q}_n+\textbf{i} \tilde{Q}_{n+1}+ \boldsymbol{\varepsilon} \tilde{Q}_{n+2}+\textbf{h} \tilde{Q}_{n+3}
\end{align}
and
\begin{align}
     \widehat{P}_n & =\breve{P}_n+ i \breve{P}_{n+1}+ j \breve{P}_{n+2}+ k \breve{P}_{n+3}\\
   & = \tilde{P}_n + \textbf{i}\tilde{P}_{n+1} + \boldsymbol{\varepsilon} \tilde{P}_{n+2} + \textbf{h} \tilde{P}_{n+3}.
\end{align}
then, the addition and subtraction of the Fibonacci hybrid quaternions are defined by
\begin{small}
    \begin{align*}
    \widehat{Q}_n \mp \widehat{P}_n & =(\breve{Q}_n+ i \breve{Q}_{n+1}+ j \breve{Q}_{n+2}+ k \breve{Q}_{n+3}) \mp (\breve{P}_n+ i \breve{P}_{n+1}+ j \breve{P}_{n+2}+ k \breve{P}_{n+3})\\
   & = \small(\breve{Q}_n \mp \breve{P}_n)+ i (\breve{Q}_{n+1} \mp \breve{P}_{n+1}) + j (\breve{Q}_{n+2} \mp \breve{P}_{n+2})+ k (\breve{Q}_{n+3} \mp \breve{P}_{n+3}),\\  
   \widehat{Q}_n \mp \widehat{P}_n & =(\tilde{Q}_n+\textbf{i} \tilde{Q}_{n+1}+ \boldsymbol{\varepsilon} \tilde{Q}_{n+2}+\textbf{h} \tilde{Q}_{n+3}) \mp (\tilde{P}_n+\textbf{i} \tilde{P}_{n+1}+ \boldsymbol{\varepsilon} \tilde{P}_{n+2}+\textbf{h} \tilde{P}_{n+3})\\
   & = (\tilde{Q}_n \mp \tilde{P}_n +\textbf{i} (\tilde{Q}_{n+1} \mp \tilde{P}_{n+1})+ \boldsymbol{\varepsilon} (\tilde{Q}_{n+2} \mp\tilde{P}_{n+2})+\textbf{h} (\tilde{Q}_{n+3} \mp \tilde{P}_{n+3}).
    \end{align*} 
\end{small}
\end{defn}
\begin{defn}
Multiplication of the Fibonacci hybrid quaternions is defined in terms of Fibonacci hybrid numbers $(\breve{Q}_n, \breve{P}_n)$ as follows:
\begin{align*}
    \widehat{Q}_n  \widehat{P}_n & =(\breve{Q}_n+ i \breve{Q}_{n+1}+ j \breve{Q}_{n+2}+ k \breve{Q}_{n+3})(\breve{P}_n+ i \breve{P}_{n+1}+ j \breve{P}_{n+2}+ k \breve{P}_{n+3})\\
   & = (\breve{Q}_n \breve{P}_n - \breve{Q}_{n+1} \breve{P}_{n+1} -  \breve{Q}_{n+2} \breve{P}_{n+2} - \breve{Q}_{n+3} \breve{P}_{n+3}) \\ & \quad + i (\breve{Q}_n \breve{P}_{n+1} + \breve{Q}_{n+1} \breve{P}_n +  \breve{Q}_{n+2} \breve{P}_{n+3} - \breve{Q}_{n+3} \breve{P}_{n+2}) \\
   & \quad + j (\breve{Q}_n \breve{P}_{n+2} -  \breve{Q}_{n+1} \breve{P}_{n+3} + \breve{Q}_{n+2} \breve{P}_n  + \breve{Q}_{n+3} \breve{P}_{n+1}) \\
   & \quad + k ( \breve{Q}_{n} \breve{P}_{n+3}+\breve{Q}_{n+1} \breve{P}_{n+2}- \breve{Q}_{n+2} \breve{P}_{n+1}+ \breve{Q}_{n+3} \breve{P}_{n})
\end{align*}
or in terms of Fibonacci quaternions $(\tilde{Q}_n,\tilde{P}_n)$ it can be defined as follows:
\begin{align*}
    \widehat{Q}_n \widehat{P}_n & =(\tilde{Q}_n+\textbf{i} \tilde{Q}_{n+1}+ \boldsymbol{\varepsilon} \tilde{Q}_{n+2}+\textbf{h} \tilde{Q}_{n+3}) (\tilde{P}_n+\textbf{i} \tilde{P}_{n+1}+ \boldsymbol{\varepsilon} \tilde{P}_{n+2}+\textbf{h} \tilde{P}_{n+3}) \\
   & = (\tilde{Q}_n \tilde{P}_n - \tilde{Q}_{n+1} \tilde{P}_{n+1} + \tilde{Q}_{n+3} \tilde{P}_{n+3} + \tilde{Q}_{n+1} \tilde{P}_{n+2} + \tilde{Q}_{n+2} \tilde{P}_{n+1}) \\
   & \quad +\textbf{i} (\tilde{Q}_n \tilde{P}_{n+1} + \tilde{Q}_{n+1} \tilde{P}_{n} + \tilde{Q}_{n+1} \tilde{P}_{n+3} - \tilde{Q}_{n+3} \tilde{P}_{n+1}) \\ 
   & \quad +\boldsymbol{\varepsilon} 
         (\tilde{Q}_n \tilde{P}_{n+2} + \tilde{Q}_{n+1} \tilde{P}_{n+3} +  \tilde{Q}_{n+2} \tilde{P}_{n}- \tilde{Q}_{n+2}  \tilde{P}_{n+3} \\ & \qquad- \tilde{Q}_{n+3}  \tilde{P}_{n+1} + \tilde{Q}_{n+3}  \tilde{P}_{n+2})\\ 
   & \quad +\textbf{h}(\tilde{Q}_n \tilde{P}_{n+3} - \tilde{Q}_{n+1} \tilde{P}_{n+2} + \tilde{Q}_{n+2} \tilde{P}_{n+1} + \tilde{Q}_{n+3} \tilde{P}_{n}).
\end{align*}
\end{defn}

The scalar and vector parts of $\widehat{Q}_n$ which is the $nth$ term of the Fibonacci hybrid quaternion sequence $(\widehat{Q}_n)$ are denoted by
\begin{equation*}
    S_{\widehat{Q}_n} = \breve{Q}_n \quad \text{and} \quad V_{\widehat{Q}_n} = i \breve{Q}_{n+1} + j \breve{Q}_{n+2} + k \breve{Q}_{n+3}
\end{equation*}
So, any Fibonacci hybrid quaternion $\widehat{Q}_n$ can be written as $\widehat{Q}_n = S_{\widehat{Q}_n} + V_{\widehat{Q}_n}$. Now we can redefine addition and subtraction as
 \begin{align*}
     \widehat{Q}_n \mp \widehat{P}_n & = (S_{\widehat{Q}_n} + V_{\widehat{Q}_n}) \mp (S_{\widehat{P}_n} + V_{\widehat{P}_n})\\
   & = (S_{\widehat{Q}_n} \mp S_{\widehat{P}_n}) + (V_{\widehat{Q}_n} \mp V_{\widehat{P}_n}),
 \end{align*}
 and multiplication
\begin{align*}
    \widehat{Q}_n \widehat{P}_n & = (S_{\widehat{Q}_n} + V_{\widehat{Q}_n}) (S_{\widehat{P}_n} + V_{\widehat{P}_n})\\
   & = S_{\widehat{Q}_n} S_{\widehat{P}_n} - \langle V_{\widehat{Q}_n}, V_{\widehat{P}_n} \rangle + S_{\widehat{Q}_n} V_{\widehat{P}_n} + S_{\widehat{P}_n} V_{\widehat{Q}_n} + V_{\widehat{Q}_n} \times V_{\widehat{P}_n}.
\end{align*}

\begin{defn}
The conjugate of Fibonacci hybrid quaternion can be define three different types for $\widehat{F}_n=\tilde{F}_n + \textbf{i} \tilde{F}_{n+1} + \boldsymbol{\varepsilon} \tilde{F}_{n+2} + \textbf{h} \tilde{F}_{n+3}$
\begin{itemize}
    \item[\textbf{i)}] Quaternion conjugate, $\overline{\widehat{F}_n}$: 
        \qquad  $
           \overline{\widehat{F}_n} = \overline{\tilde{F}}_n+\textbf{i} \overline{\tilde{F}}_{n+1}+\boldsymbol{\varepsilon} \overline{\tilde{F}}_{n+2}+\textbf{h} \overline{\tilde{F}}_{n+3}
         $,
    \item[\textbf{ii)}] Hybrid conjugate, $(\widehat{F}_n)^C$: \qquad $
    (\widehat{F}_n)^C=\tilde{F}_n-\textbf{i} \tilde{F}_{n+1} - \boldsymbol{\varepsilon}\tilde{F}_{n+2} - \textbf{h} \tilde{F}_{n+3}    
    $,
    \item[\textbf{iii)}] Total conjugate, $(\widehat{F}_n)^{\dag}$: 
    \qquad $(\widehat{F}_n)^{\dag}=\overline{(\widehat{F}_n)^C}=\overline{\tilde{F}}_n- \textbf{i} \overline{\tilde{F}}_{n+1}- \boldsymbol{\varepsilon} \overline{\tilde{F}}_{n+2}-\textbf{h} \overline{\tilde{F}}_{n+3}$.
\end{itemize}
\end{defn}

\begin{thm}
Let $\widehat{F}_n$ be $n$th term of the Fibonacci sequence. Then, for $n\geq1$ we can give the following relations:
\begin{itemize}
    \item[\textbf{i)}] $\widehat{F}_n+\widehat{F}_{n+1}= \widehat{F}_{n+2}$,
    \item[\textbf{ii)}] $\widehat{F}_n-i\widehat{F}_{n+1}-j\widehat{F}_{n+2}-k\widehat{F}_{n+3}=\breve{F}_n+\breve{F}_{n+2}+\breve{F}_{n+4}+\breve{F}_{n+6}=\breve{L}_{n+1}+\breve{L}_{n+5}$,
    \item[\textbf{iii)}] $\widehat{F}_n-\textbf{i} \widehat{F}_{n+1}-\boldsymbol{\varepsilon} \widehat{F}_{n+2}-\textbf{h} \widehat{F}_{n+3}=\tilde{F}_n-\tilde{F}_{n+2}-2\tilde{F}_{n+3}+\tilde{F}_{n+6}$.
\end{itemize}
\end{thm}

\textbf{\textit{Proof.}}
Let $\breve{F}_n,$ $\breve{L}_n$ and $\tilde{F}_n$ be $n$th Fibonacci hybrid number, $n$th Lucas hybrid number and $n$th Fibonacci quaternion, respectively.

    \noindent \textbf{i)} We can show this equality in two ways; first one is using Fibonacci hybrid numbers: 
        \begin{align*}
           \widehat{F}_n + \widehat{F}_{n+1} & = (\breve{F}_n+ i \breve{F}_{n+1}+ j \breve{F}_{n+2}+ k \breve{F}_{n+3}) +(\breve{F}_{n+1}+ i \breve{F}_{n+2}+ j \breve{F}_{n+3}+ k \breve{F}_{n+4})\\
           & = (\breve{F}_n+\breve{F}_{n+1}) + i (\breve{F}_{n+1}+\breve{F}_{n+2}) + j (\breve{F}_{n+2}+\breve{F}_{n+3}) + k (\breve{F}_{n+3}+\breve{F}_{n+4}) \\ 
           & = \breve{F}_{n+2} + i \breve{F}_{n+3} + j \breve{F}_{n+4} + k \breve{F}_{n+5}= \widehat{F}_{n+2} .
        \end{align*}
The second one is using Fibonacci quaternions:
 \begin{align*}
   \widehat{F}_n + \widehat{F}_{n+1} & = (\tilde{F}_n+\textbf{i} \tilde{F}_{n+1}+\boldsymbol{\varepsilon} \tilde{F}_{n+2}+\textbf{h} \tilde{F}_{n+3})
    +(\tilde{F}_{n+1}+\textbf{i} \tilde{F}_{n+2}+\boldsymbol{\varepsilon} \tilde{F}_{n+3}+\textbf{h} \tilde{F}_{n+4})\\
   & = (\tilde{F}_n+\tilde{F}_{n+1})+\textbf{i} (\tilde{F}_{n+1}+\tilde{F}_{n+2}
    + \boldsymbol{\varepsilon} (\tilde{F}_{n+2}+\tilde{F}_{n+3})+ \textbf{h} (\tilde{F}_{n+3}+\tilde{F}_{n+4}) \\
   & = \tilde{F}_{n+2}+\textbf{i} \tilde{F}_{n+3}+ \boldsymbol{\varepsilon} \tilde{F}_{n+4}+ \textbf{h} \tilde{F}_{n+5} = {\widehat{F}_{n+2}}.
\end{align*}

   \noindent \textbf{ii)} If we use $\breve{F}_{n-1}+\breve{F}_{n+1}=\breve{L}_n$ \cite{horhyb}, we have
    \begin{align*}
       \widehat{F}_n - i \widehat{F}_{n+1} - j \widehat{F}_{n+2}- k \widehat{F}_{n+3} & = (\breve{F}_n+i\breve{F}_{n+1}+j\breve{F}_{n+2}+k\breve{F}_{n+3})-i(\breve{F}_{n+1}+i\breve{F}_{n+2}+j\breve{F}_{n+3}+k\breve{F}_{n+4}) \\
       &\quad -j(\breve{F}_{n+2}+i\breve{F}_{n+3}+j\breve{F}_{n+4}+k\breve{F}_{n+5}) -k(\breve{F}_{n+3}+i\breve{F}_{n+4}+j\breve{F}_{n+5}+k\breve{F}_{n+6}) \\
       & =\breve{F}_n+\breve{F}_{n+2}+\breve{F}_{n+4}+\breve{F}_{n+6} \\
       & =\breve{L}_{n+1}+\breve{L}_{n+5}.
\end{align*}
\begin{align*}
\noindent \textbf{iii)} \widehat{F}_n-\textbf{i} \widehat{F}_{n+1}-\boldsymbol{\varepsilon} \widehat{F}_{n+2}-h \widehat{F}_{n+3}& = (\widehat{F}_n+ \textbf{i} \widehat{F}_{n+1}+\boldsymbol{\varepsilon} \widehat{F}_{n+2}+\textbf{h}\widehat{F}_{n+3}) -\textbf{i} (\widehat{F}_{n+1}+\textbf{i}\widehat{F}_{n+2}+\boldsymbol{\varepsilon} \widehat{F}_{n+3}+\textbf{h}\widehat{F}_{n+4}) \\
   & -\boldsymbol{\varepsilon}(\widehat{F}_{n+2}+\textbf{i}\widehat{F}_{n+3}+\boldsymbol{\varepsilon} \widehat{F}_{n+4}+\textbf{h}\widehat{F}_{n+5}) -\textbf{h}(\widehat{F}_{n+3}+\textbf{i}\widehat{F}_{n+4}+\boldsymbol{\varepsilon} \widehat{F}_{n+5}+\textbf{h}\widehat{F}_{n+6}) \\
   & =\widehat{F}_n-\widehat{F}_{n+2}-2\widehat{F}_{n+3}+\widehat{F}_{n+6}.
\end{align*}

\begin{thm}
Let $\widehat{F}_n$ and $\widehat{L}_n$ be $nth$ Fibonacci hybrid quaternion and the Lucas hybrid quaternion sequences respectively. The following relations are satisfied:
\begin{itemize}
    \item[\textbf{i)}] $\widehat{F}_{n-1}+\widehat{F}_{n+1}=\widehat{L}_n$
    \item[\textbf{ii)}] $\widehat{F}_{n+2}-\widehat{F}_{n-2}=\widehat{L}_n$
\end{itemize}
\end{thm}
\textbf{\textit{Proof.}}
\begin{itemize}
    \item[\textbf{i)}] We know that $\breve{L}_n=\breve{F}_{n-1}+\breve{F}_{n+1}$. Then we have
\begin{align*}
   \widehat{F}_{n-1}+ \widehat{F}_{n+1} & = (\breve{F}_{n-1}+ i \breve{F}_{n}+j \breve{F}_{n+1}+k \breve{F}_{n+2}) +(\breve{F}_{n+1}+ i \breve{F}_{n+2}+j \breve{F}_{n+3}+k \breve{F}_{n+4})\\
   & = (\breve{F}_{n-1}+\breve{F}_{n+1})+ i (\breve{F}_{n}+F_{n+2}) +j ( \breve{F}_{n+1}+\breve{F}_{n+3})+k (\breve{F}_{n+2}+\breve{F}_{n+4}) \\
   & = \breve{L}_n+ i \breve{L}_{n+1}+j \breve{L}_{n+2}+k \breve{L}_{n+3} \\
   & =\widehat{L}_n
\end{align*}
    \item[\textbf{ii)}] We know that $\breve{F}_n=\breve{F}_{n+2}-\breve{F}_{n-2}$.
\begin{align*}
   \widehat{F}_{n+2}- \widehat{F}_{n-2} & = (\breve{F}_{n+2}+ i \breve{F}_{n+3}+j \breve{F}_{n+4}+k \breve{F}_{n+5}) - (\breve{F}_{n-2}+ i \breve{F}_{n-1}+j \breve{F}_{n}+k \breve{F}_{n+1})\\
   & = (\breve{F}_{n+2}-\breve{F}_{n-2})+ i (\breve{F}_{n+3}-\breve{F}_{n-1}) +j ( \breve{F}_{n+4}-\breve{F}_{n})+k (\breve{F}_{n+5}-\breve{F}_{n+1}) \\
   & = \breve{L}_n+ i \breve{L}_{n+1}+j \breve{L}_{n+2}+k \breve{L}_{n+3} \\
   & =\widehat{L}_n
\end{align*}
\end{itemize}

\begin{thm}
Let $\widehat{F}_n$ be $nth$ Fibonacci hybrid quaternion sequence ($\widehat{F}_n$) and $\overline{\widehat{F}_n}$, $\widehat{F}_n^C$, $\widehat{F}_n^{\dag}$ be the quaternion conjugate, hybrid conjugate and total conjugate of $nth$ Fibonacci hybrid quaternion respectively. The following relations are satisfies:
\begin{itemize}
    \item[\textbf{i)}] $\widehat{ F}_n+\overline{\widehat{ F}_n}=2 \breve{F}_n$
    \item[\textbf{ii)}] $\widehat{ F}_n+\widehat{ F}_n^C=2\tilde{F}_n$
    \item[\textbf{iii)}] $\widehat{ F}_n+\widehat{ F}_n^{\dag}= -2 F_n - 8 F_{n+1}+2(\widehat{ F}_{n+1}+\widehat{ F}_{n+2}+\widehat{ F}_{n+3})$
\end{itemize}
\end{thm}

\textbf{\textit{Proof.}}
\begin{align*}
 \textbf{i)} \qquad \widehat{ F}_n+\overline{\widehat{ F}_n} & = (\breve{F}_{n}+ i \breve{F}_{n+1}+j \breve{F}_{n+2}+k \breve{F}_{n+3})+(\breve{F}_{n}- i \breve{F}_{n+1}-j \breve{F}_{n+2}-k \breve{F}_{n+3})\qquad \qquad \qquad\\
   & = 2\breve{F}_n
\end{align*}
\begin{align*}
  \textbf{ii)} \qquad  \widehat{ F}_n+\widehat{ F}_n^C & = (\breve{F}_{n}+ i \breve{F}_{n+1}+j \breve{F}_{n+2}+k \breve{F}_{n+3}) +( {\breve{F}_{n}^C}+ i  {\breve{F}_{n+1}^C}+j  {\breve{F}_{n+2}^C}+k  {\breve{F}_{n+3}^C})\qquad \qquad \qquad \\
   & = (\breve{F}_{n}+ {\breve{F}_{n}^C})+ i (\breve{F}_{n+1}+ {\breve{F}_{n+1}^C}) + j (\breve{F}_{n+2}+ {\breve{F}_{n+2}^C})+  k (\breve{F}_{n+3}+ {\breve{F}_{n+3}^C})\\
   & = 2(F_n+i F_{n+1}+ j F_{n+2} k F_{n+3})\\
   & = 2 \tilde{F}_n
\end{align*}
\begin{align*}
  \textbf{iii)} \qquad   \widehat{ F}_n+\widehat{ F}_n^{\dag} & = ( \breve{F}_{n}+ i  \breve{F}_{n+1}+j  \breve{F}_{n+2}+k  \breve{F}_{n+3})+ ( { \breve{F}_{n}^C}- i  { \breve{F}_{n+1}^C}-j  { \breve{F}_{n+2}^C}-k  { \breve{F}_{n+3}^C})\qquad \qquad \qquad \\ 
   & = ( \breve{F}_{n}+ { \breve{F}_{n}^C})+ i ( \breve{F}_{n+1}- { \breve{F}_{n+1}^C}) + j ( \breve{F}_{n+2}- { \breve{F}_{n+2}^C})+  k ( \breve{F}_{n+3}- { \breve{F}_{n+3}^C})\\
   & = -2F_n-8F_{n+1}+2( \breve{ F}_{n+1}+ \breve{ F}_{n+2}+ \breve{ F}_{n+3}) \qquad \qquad \qquad
\end{align*}

\begin{thm}[Binet's Formulas]
Let $\widehat{F}_n$ and $\widehat{L}_n$ be Fibonacci hybrid quaternion and Lucas hybrid quaternion respectively. The Binet formulas for these hybrid quaternions are given as follows:
\smallskip
\begin{itemize}
    \item[\textbf{i)}] $\widehat{F}_n=\frac{\alpha^*\underline{\alpha}\alpha^n-\beta^*\underline{\beta}\beta^n}{\alpha-\beta}$
    \smallskip
    \item[\textbf{ii)}]
    $\widehat{L}_n=\alpha^*\underline{\alpha}\alpha^n+\beta^*\underline{\beta}\beta^n$
\end{itemize}
\smallskip
where  $\alpha^*=1+\textbf{i}\alpha+\boldsymbol{\varepsilon}\alpha^2+\textbf{h}\alpha^3$, $\beta^*=1+\textbf{i}\beta+\boldsymbol{\varepsilon}\beta^2+\textbf{h}\beta^3$, $\underline{\alpha}=1+i\alpha+j\alpha^2+k\alpha^3$ and $\underline{\beta}=1+i\beta+j\beta^2+k\beta^3$.
\end{thm}

\textbf{\textit{Proof.}}
In \cite{fibquat}, Halici gave the Binet's formula for Fibonacci and Lucas quaternions by
    \begin{equation}\label{23}
        \tilde{F}_n=\frac{\underline{\alpha}\alpha^n-\underline{\beta}\beta^n}{\alpha-\beta}
    \end{equation}
and
    \begin{equation}\label{24}
        \tilde{L}_n=\underline{\alpha}\alpha^n+\underline{\beta}\beta^n.
    \end{equation}
\textbf{i)} By using \eqref{23}, we have
\begin{align*}
   \widehat{F}_n & =  \tilde{F}_n+\textbf{i} \tilde{F}_{n+1}+\boldsymbol{\varepsilon}  \tilde{F}_{n+2}+\textbf{h}  \tilde{F}_{n+3}\\
   & = \frac{\underline{\alpha}\alpha^n-\underline{\beta}\beta^n}{\alpha-\beta}+\textbf{i} \frac{\underline{\alpha}\alpha^{n+1}-\underline{\beta}\beta^{n+1}}{\alpha-\beta}+\boldsymbol{\varepsilon} \frac{\underline{\alpha}\alpha^{n+2}-\underline{\beta}\beta^{n+2}}{\alpha-\beta}+\textbf{h} \frac{\underline{\alpha}\alpha^{n+3}-\underline{\beta}\beta^{n+3}}{\alpha-\beta} \\
   & = \frac{(\underline{\alpha}\alpha^n-\underline{\beta}\beta^n)+\textbf{i}(\underline{\alpha}\alpha^{n+1}-\underline{\beta}\beta^{n+1})+\boldsymbol{\varepsilon}(\underline{\alpha}\alpha^{n+2}-\underline{\beta}\beta^{n+2})+\textbf{h}(\underline{\alpha}\alpha^{n+3}-\underline{\beta}\beta^{n+3})}{\alpha-\beta} \\
   & = \frac{(\underline{\alpha}\alpha^n+\textbf{i}\underline{\alpha}\alpha^{n+1}+\boldsymbol{\varepsilon} \underline{\alpha}\alpha^{n+2}+\textbf{h} \underline{\alpha}\alpha^{n+3})-(\underline{\beta}\beta^n+\textbf{i}\underline{\beta}\beta^{n+1}+\boldsymbol{\varepsilon} \underline{\beta}\beta^{n+2}+\textbf{h} \underline{\beta}\beta^{n+3})}{\alpha-\beta} \\
   & = \frac{\underline{\alpha}\alpha^n(1+\textbf{i}\alpha+\boldsymbol{\varepsilon}\alpha^2+\textbf{h}\alpha^3)-\underline{\beta}\beta^n(1+\textbf{i}\beta+\boldsymbol{\varepsilon}\beta^2+\textbf{h}\beta^3)}{\alpha-\beta} \\
   & = \frac{\alpha^*\underline{\alpha}\alpha^n-\beta^*\underline{\alpha}\beta^n}{\alpha-\beta}. 
   \end{align*}

\noindent \textbf{ii)} By using \eqref{24}, we have
\begin{align*}
   \widehat{L}_n & = \tilde{L}_n+\textbf{i}\tilde{L}_{n+1}+\boldsymbol{\varepsilon} \tilde{L}_{n+2}+h \tilde{L}_{n+3}\\
   & = \underline{\alpha}\alpha^n+\underline{\beta}\beta^n+\textbf{i} \underline{\alpha}\alpha^{n+1}+\underline{\beta}\beta^{n+1}+\boldsymbol{\varepsilon} \underline{\alpha}\alpha^{n+2} +\underline{\beta}\beta^{n+2}+\textbf{h} \underline{\alpha}\alpha^{n+3}+\underline{\beta}\beta^{n+3} \\
   & = (\underline{\alpha}\alpha^n+\underline{\beta}\beta^n)+\textbf{i}(\underline{\alpha}\alpha^{n+1}+\underline{\beta}\beta^{n+1}) +\boldsymbol{\varepsilon}(\underline{\alpha}\alpha^{n+2}+\underline{\beta}\beta^{n+2})+\textbf{h}(\underline{\alpha}\alpha^{n+3}+\underline{\beta}\beta^{n+3}) \\
   & = (\underline{\alpha}\alpha^n+\textbf{i}\underline{\alpha}\alpha^{n+1}+\boldsymbol{\varepsilon} \underline{\alpha}\alpha^{n+2}+\textbf{h} \underline{\alpha}\alpha^{n+3}) +(\underline{\beta}\beta^n+\textbf{i}\underline{\beta}\beta^{n+1}+\boldsymbol{\varepsilon} \underline{\beta}\beta^{n+2}+\textbf{h} \underline{\beta}\beta^{n+3}) \\
   & = \underline{\alpha}\alpha^n(1+\textbf{i}\alpha+\boldsymbol{\varepsilon}\alpha^2+\textbf{h}\alpha^3)+\underline{\beta}\beta^n(1+\textbf{i}\beta+\boldsymbol{\varepsilon}\beta^2+\textbf{h}\beta^3) \\
   & = \alpha^*\underline{\alpha}\alpha^n+\beta^*\underline{\beta}\beta^n.
\end{align*}

\begin{thm}[Cassini's Identities]
The following equations are hold:
\begin{align*}
    C_1&=\widehat{F}_{n+1}\widehat{F}_{n-1}-\widehat{F}_n^2=(-1)^n\frac{\alpha\alpha^*\beta^*\underline{\alpha}\underline{\beta}-\beta\beta^*\alpha^*\underline{\beta}\underline{\alpha}}{\alpha-\beta},\\
    C_2&=\widehat{L}_{n+1}\widehat{L}_{n-1}-\widehat{L}_n^2=(-1)^n\sqrt{5}(\alpha\alpha^*\beta^*\underline{\alpha}\underline{\beta}-\beta\beta^*\alpha^*\underline{\beta}\underline{\alpha})
\end{align*}
where $\alpha$ and $\beta$ be the roots of equation $x^2-2x-1=0$, $\underline{\alpha}=1+\alpha i+ \alpha^2 j + \alpha^3 L$, $\underline{\beta}=1+\beta i+ \beta^2 j + \beta^3 k$, $\alpha^*=1+\alpha\mathbf{i}+\alpha^2\boldsymbol{\varepsilon}+\alpha^3\mathbf{h}$ and $\beta^*=1+\beta\mathbf{i}+\beta^2\boldsymbol{\varepsilon}+\beta^3\mathbf{h}$.
\end{thm}

\textbf{\textit{Proof.}} For the first Cassini identity $C_1$, we get
\begin{align*}
    C_1 
    &= \bigg(\frac{\alpha^{n+1}\alpha^*\underline{\alpha}-\beta^{n+1}\beta^*\underline{\beta}}{\alpha-\beta}\bigg)\bigg(\frac{\alpha^{n-1}\alpha^*\underline{\alpha}-\beta^{n-1}\beta^*\underline{\beta}}{\alpha-\beta} \bigg)- \bigg(\frac{\alpha^{n}\alpha^*\underline{\alpha}-\beta^{n}\beta^*\underline{\beta}}{\alpha-\beta}\bigg)^2 \\
    &= \frac{\alpha^{2n}(\alpha^*)^2\underline{\alpha}^2-\alpha^{n+1}\beta^{n-1}\alpha^*\beta^*\underline{\alpha}\underline{\beta}-\beta^{n+1}\alpha^{n-1}\beta^*\alpha^*\underline{\beta}\underline{\alpha}+\beta^{2n}(\beta^*)^2\beta^2}{(\alpha-\beta)^2} \\
    & \quad - \frac{\alpha^{2n}(\alpha^*)^2\underline{\alpha}^2-\alpha^{n}\beta^{n}\alpha^*\beta^*\underline{\alpha}\underline{\beta}-\beta^{n}\alpha^{n}\beta^*\alpha^*\underline{\beta}\underline{\alpha}+\beta^{2n}(\beta^*)^2\beta^2}{(\alpha-\beta)^2} \\
    &= \frac{\alpha^{n-1}\beta^n\beta^*\alpha^*(\alpha-\beta)\underline{\beta}\underline{\alpha}-\alpha^n\beta^{n-1}\alpha^*\beta^*(\alpha-\beta)\underline{\alpha}\underline{\beta}}{(\alpha-\beta)^2} \\
    &= (-1)^n\frac{(\alpha\alpha^*\beta^*\underline{\alpha}\underline{\beta}-\beta\beta^*\alpha^*\underline{\beta}\underline{\alpha})}{\alpha-\beta}
\end{align*}
For the second Cassini identity $C_2$, we get
\begin{align*}
    C_2 &= (\alpha^{n+1}\alpha^*\underline{\alpha}+\beta^{n+1}\beta^*\underline{\beta})(\alpha^{n-1}\alpha^*\underline{\alpha}+\beta^{n-1}\beta^*\underline{\beta})  - (\alpha^{n}\alpha^*\underline{\alpha}+\beta^{n}\beta^*\underline{\beta})^2 \\
    &= \alpha^{2n}(\alpha^*)^2\underline{\alpha}^2+\alpha^{n+1}\beta^{n-1}\alpha^*\beta^*\underline{\alpha}\underline{\beta}+\beta^{n+1}\alpha^{n-1}\beta^*\alpha^*\underline{\beta}\underline{\alpha}+\beta^{2n}(\beta^*)^2\beta^2 \\
    & \quad - \alpha^{2n}(\alpha^*)^2\underline{\alpha}^2+\alpha^{n}\beta^{n}\alpha^*\beta^*\underline{\alpha}\underline{\beta}+\beta^{n}\alpha^{n}\beta^*\alpha^*\underline{\beta}\underline{\alpha}+\beta^{2n}(\beta^*)^2\beta^2 \\
    &= \alpha^{n-1}\beta^n\beta^*\alpha^*(\beta-\alpha)\underline{\beta}\underline{\alpha}+\alpha^n\beta^{n-1}\alpha^*\beta^*(\alpha-\beta)\underline{\alpha}\underline{\beta} \\
    &= \alpha^{n-1}\beta^{n-1}(\alpha-\beta)(\alpha\alpha^*\beta^*\underline{\alpha}\underline{\beta}-\beta\beta^*\alpha^*\underline{\beta}\underline{\alpha}) \\
    &= (-1)^n\sqrt{5}(\alpha\alpha^*\beta^*\underline{\alpha}\underline{\beta}-\beta\beta^*\alpha^*\underline{\beta}\underline{\alpha}).
\end{align*}

\section{Conclusion}
In this paper, we have introduced the Horadam hybrid quaternions and some special classes of number sequences such as Fibonacci, Lucas, Pell and Jacobsthal hybrid quaternions. Especially, we have examined  the Fibonacci and Lucas hybrid quaternions comprehensively. Moreover, some identities like Binet's formula and Cassini's identity for Fibonacci and Lucas hybrid quaternions have been given. Other sequences on Hybrid quaternions can be studied as future works.

\end{document}